\newif\ifsmfart
\IfFileExists{smfart.cls}
   {\documentclass[12pt,english]{smfart}
    \IfFileExists{smfenum.sty}{\usepackage{smfenum}}{}
    \usepackage{bull}
    \smfarttrue}
{\message{^^J*** It would be better to typeset this file with smfart.cls ***^^J^^J}
    
\documentclass[12pt]{amsart}}

\usepackage{amssymb}
\usepackage{epsfig}
\usepackage{graphicx}
\numberwithin{equation}{section}

\input xy
\xyoption{all}

\advance\headheight 2pt
\calclayout
\allowdisplaybreaks[3]

\theoremstyle{plain}
\newtheorem{prop}[subsection]{Proposition}

\newtheorem{theo}[subsection]{Theorem}
\newtheorem{coro}[subsection]{Corollary}

\newtheorem{lem}[subsection]{Lemma}

\newtheorem{defn}[subsection]{Definition}

\theoremstyle{definition}

\newtheorem{conj}[subsection]{Conjecture}

\theoremstyle{remark}
\newtheorem{rem}[subsection]{Remark}

\newtheorem{exam}[subsection]{Example}

\newcommand{\bA}{{\mathbb A}}

\newcommand{\bG}{{\mathbb G}}
\newcommand{\bN}{{\mathbb N}}

\newcommand{\bR}{{\mathbb R}}
\newcommand{\bP}{{\mathbb P}}
\newcommand{\bZ}{{\mathbb Z}}
\newcommand{\cF}{{\mathcal F}}
\newcommand{\cO}{{\mathcal O}}
\def\cL{{\mathcal L}}
\def\cX{{\mathfrak X}}
\newcommand{\ra}{{\rightarrow}}
\newcommand{\hra}{{\hookrightarrow}}
\newcommand{\Spec}{{\rm Spec}}
\newcommand{\Move}{{\rm Mov}}
\newcommand{\Proj}{{\rm Proj}}

\newcommand{\tS}{\tilde S}
\newcommand{\NE}{\mathrm{NE}}
\newcommand{\NM}{\mathrm{NM}}
\newcommand{\Pic}{\mathrm{Pic}}
\newcommand{\Hom}{\mathrm{Hom}}
\newcommand{\Cone}{\mathrm{Cone}}
\newcommand{\Cox}{\mathrm{Cox}}

\newcommand{\Dfour}{{\mathbf D}_4}
\newcommand{\Esix}{{\mathbf E}_6}
\def\la{{\lambda}}
\def\rA{{\mathrm A}}
\def\rN{{\mathrm N}}
\def\Tr{{\rm Tr}}
\def\Gal{{\rm Gal}}

\makeatother
\makeatletter

\author{Brendan Hassett}
\address{Department of Mathematics \\
Rice University \\
Houston, TX 77005}
\email{hassett@rice.edu}
\author{Yuri Tschinkel}
\address{Department of Mathematics \\
         Princeton University\\ 
         Fine Hall, Washington Road\\
         Princeton, NJ 08544-1000,  U.S.A.}
\email{ytschink@math.princeton.edu}
\thanks{The first author was 
partially supported by the Sloan Foundation and by
NSF Grants 0196187 and 0134259. 
The second author was partially supported by NSF Grant 0100277.}

\usepackage{graphicx}

\title[Universal torsors and Cox rings]{Universal torsors and Cox rings}
\date{\today}
\begin{document}

\begin{abstract}
We study the equations of universal torsors on rational surfaces. 
\end{abstract}

\maketitle
\tableofcontents

\setcounter{section}{0}
\section*{Introduction}
\label{sect:introduction}

The study of surfaces over nonclosed fields $k$ leads naturally
to certain auxiliary varieties, called {\em universal torsors}.
The proofs of the Hasse principle and weak approximation
for certain Del Pezzo surfaces
required a very detailed knowledge of the projective geometry, in fact,
explicit equations, for these torsors \cite{CTS}, \cite{CSS}, \cite{CSS2},
\cite{SS}, \cite{S}, \cite{Sb}.
More recently, Salberger proposed using universal torsors to count
rational points of bounded
height, obtaining the first sharp upper bounds on split Del Pezzo surfaces
of degree 5 and asymptotics on split toric varieties over $\mathbb Q$
\cite{S-ast}. This approach was further developed in the work of
Peyre, de la Bret\`eche, and Heath-Brown 
\cite{Peyre}, \cite{Peyre2}, \cite{Br}, \cite{HB-C}.

\

Colliot-Th\'el\`ene and Sansuc have given a general formalism for
writing down equations for these torsors.
We briefly sketch their method:
Let $X$ be a smooth projective variety and $\{ D_j\}_{j\in J} $ a finite set
of irreducible divisors on $X$ such that
$U:=X\setminus  \cup_{j\in J} D_j$ has trivial Picard group.
In practice, one usually chooses generators
of the effective cone of $X$, e.g., the lines on the Del Pezzo surface.
Consider the resulting exact sequence:
$$
0\longrightarrow \bar{k}[U]^*/\bar{k}^*\longrightarrow \oplus_{j\in J}
\bZ D_j \longrightarrow \Pic(X_{\bar{k}})\longrightarrow 0.
$$
Applying $\Hom(-,\bG_m)$, one obtains
an exact sequence of tori
$$
1\longrightarrow T(X)\longrightarrow T\longrightarrow R\longrightarrow 1,
$$
where the first term is the {\em N\'eron-Severi torus} of $X$.
Suppose we have a collection
of rational functions, invertible on $U$, which form a basis for
the relations among the $\{D_j\}_{j\in J}$.
These can be interpreted as a section
$U\ra R\times U$, and thus naturally
induce a $T(X)$-torsor over $U$, which  canonically extends
to the universal torsor over $X$.  In practice,
this extension can be made explicit,
yielding equations for the universal torsor.

However, when the cone generated by $\{D_j\}_{j\in J}$ is simplicial, there
are {\em no} relations and this method
gives little information.  In this paper, we outline an alternative
approach to the construction of universal torsors
and illustrate it in specific examples
where the effective cone of $X$ is simplicial.

\

We will work with varieties $X$ such that 
the Picard and the N\'eron-Severi groups
of $X$ coincide and such that the ring
$$
\Cox(X):=\bigoplus_{L\in \Pic(X)} \Gamma(X,L),
$$
is finitely generated.  This ring admits a natural action of 
the N\'eron-Severi torus and the corresponding affine variety
is a natural embedding of the universal torsor of $X$.
The challenge is to 
actually compute $\Cox(X)$ in specific examples;
Cox has shown that it is
a polynomial ring precisely when $X$ is toric \cite{Cox}.
  
\
 
Here is a roadmap of the paper: In Section~\ref{sect:general} 
we introduce Cox rings and discuss their general properties.
Finding generators for the Cox ring entails embedding
the universal torsor into affine space, which yields embeddings
of our original variety into toric quotients of this affine space.
We have collected several useful facts about toric varieties in 
Section~\ref{sect:toric}.  Section~\ref{sect:e6} is devoted to 
a detailed analysis of the unique cubic surface $S$ with an isolated 
singularity of type $\Esix$.  We compute the (simplicial) effective
cone of its minimal desingularization $\tilde{S}$, 
and produce 10 distinguished sections in $\Cox(\tilde{S})$.
These satisfy a unique equation and we show the 
universal torsor naturally embeds in 
the corresponding hypersurface in ${\mathbb A}^{10}$. 
More precisely, we get a homomorphism from the coordinate ring
of $\bA^{10}$ to $\Cox(\tilde{S})$ and the main point 
is to prove its surjectivity.   Here we use
an embedding of $\tilde{S}$ into a simplicial toric threefold $Y$, 
a quotient of ${\mathbb A}^{10}$ 
under the action of the N\'eron-Severi torus
so that $\Cox(Y)$ is the polynomial ring over the above 10 generators. 
The induced restriction map on the level of Picard groups
is an isomorphism respecting the moving cones.  We conclude
surjectivity for each degree by finding an appropriate 
birational projective model of $Y$ and using 
vanishing results on it. 
Finally, in Section~\ref{sect:d4} we write down equations 
for the universal torsors (the Cox rings) of
a split and a nonsplit cubic surface 
with an isolated singularity of 
type $\Dfour$. 

For an account of the Cox rings of smooth Del Pezzo surfaces,
we refer the reader to the paper of Batyrev and Popov 
\cite{BP}.

\

\noindent {\bf Acknowledgments:}
The results of this paper have been reported at the
American Institute of Mathematics conference ``Rational and integral
points on higher dimensional varieties''.  We benefited from the
comments of the other participants, in particular, V. Batyrev and
J.L. Colliot-Th\'el\`ene.  We also thank S. Keel for several helpful
discussions about Cox rings and M. Thaddeus for advice about the
geometric invariant theory of toric varieties.

\

\section{Generalities on the Cox ring}
\label{sect:general}

For any finite subset $\Xi$ of a real vector
space, let $\Cone(\Xi)$ denote the closed
cone generated by $\Xi$.  

Let $X$ be a normal projective variety of dimension $n$ over an
algebraically closed field $k$ of characteristic zero.  Let 
$\rA_{n-1}(X)$ and $\rN_{n-1}(X)$ denote Weil divisors on $X$ up to 
linear and numerical equivalence, respectively.  
Let $\rA_1(X)$ and $\rN_1(X)$ denote the classes of curves
up to equivalence.  
Let $\NE_{n-1}(X) \subset \rN_{n-1}(X)_{\bR}$ denote the
cone of (pseudo)effective divisors, i.e., the 
smallest real closed cone containing all the
effective divisors of $X$.  
Let $\NE_1(X)\subset \rN_1(X)_{\bR}$ denote
the cone of effective curves and   
$\NM^1(X) \subset \rN_{n-1}(X)_{\bR}$ 
the cone of nef Cartier divisors, which is
dual to the cone of effective Cartier divisors.  
By Kleiman's criterion,
this is the smallest real closed cone containing all
ample divisors of $X$.  

\

Let $L_1,\ldots,L_r$ be invertible sheaves on $X$.
For $\nu=(n_1,\ldots,n_r)\in \bN^r$ write
$$L^{\nu}:=L_1^{\otimes n_1} \otimes \ldots \otimes L_r^{\otimes n_r}.$$
Consider the ring
$$R(X,L_1,\ldots,L_r):=\bigoplus_{\nu\in \bN^r}\Gamma(X,L^{\nu}),$$
which need not be finitely generated in general.

By definition, an invertible sheaf $L$ on $X$ is {\em semiample}
if $L^N$ is globally generated for some $N>0$:

\begin{prop}{(\cite{HK}, Lemma 2.8)}
\label{prop:coxfg} 
If $L_1,\ldots,L_r$ are semiample then
$R(X,L_1,\ldots,L_r)$ is finitely generated.
\end{prop}

\begin{rem}
If the $L_i$ are ample then, after replacing each
$L_i$ by a large multiple,  $R(X,L_1,\ldots,L_r)$ is
generated by
$$
\Gamma(X,L_1) \otimes \ldots \otimes \Gamma(X,L_r).
$$ 
However, this is not generally the case if the $L_i$
are only semiample (despite the assertion in
the second part of Lemma 2.8 of \cite{HK}).  
Indeed, let $X\ra \bP^1\times \bP^1$ 
be a double cover and $L_1$ and $L_2$ be the pull-backs
of the polarizations on the $\bP^1$'s to $X$.  For suitably
large $n_1$ and $n_2$, $L_1^{n_1}\otimes L_2^{n_2}$ is very
ample and its sections embed $X$.  However, 
$$
\Gamma(X,L^{n_1}_1) \otimes \Gamma(X,L^{n_2}_2) \simeq
\Gamma(\bP^1,\cO_{\bP^1}(n_1)) \otimes \Gamma(\bP^1,\cO_{\bP^1}(n_2)),
$$
and any morphism induced by these sections factors through 
$\bP^1\times \bP^1$.   
\end{rem}

\begin{prop}
\label{prop:hu} 
Let $L_1,\ldots, L_r$ be a set of invertible sheaves on $X$
such that $L_j$ is generated by sections $s_{j,0},...,s_{j,d_j}$. 
Assume that the induced morphism $X\ra \prod_{j} \bP^{d_j}$ 
is birational into its image. 
Then the ring generated by the $s_{j,k}$'s 
has the same fraction field as $R(X,L_1,\ldots, L_r)$. 
\end{prop}

\begin{proof}
Both rings have fraction field $k(X)(t_1,...,t_r)$,
where $t_j$ is a nonzero section of $L_j$.  
\end{proof}

\

\begin{defn}
\cite{HK} 
Let $X$ be a nonsingular projective
variety so that $\Pic(X)$ is a free abelian group of rank $r$.  
The {\em Cox ring} for $X$ is defined
$$
\Cox(X):=R(X,L_1,\ldots,L_r),
$$
where $L_1,\ldots,L_r$ are lines bundles so that
\begin{enumerate}
\item{the $L_i$ form a $\bZ$-basis of $\Pic(X)$;}
\item{the cone $\Cone(\{L_1,\ldots,L_r\})$ contains $\NE_{n-1}(X)$.}
\end{enumerate}
This ring is naturally graded by $\Pic(X)$:
for $\nu\in \Pic(X)$ the $\nu$-graded piece
is denoted $\Cox(X)_{\nu}$.  
\end{defn}

\begin{prop}
\label{prop:unique} \cite{HK}
The ring $\Cox(X)$ does not depend on the choice of generators
for $\Pic(X)$.  
\end{prop}

\begin{proof} 
Consider two sets of generators
$L_1,\ldots,L_r$ and $M_1,\ldots,M_r$.  
Since $\Cone(\{L_i\})$ and $\Cone(\{M_i\})$ contain
all the effective divisors, the nonzero graded pieces
of both $R(X,L_1,\ldots,L_r)$ and
$R(X,M_1,\ldots,M_r)$ are indexed by the effective divisor
classes in $\Pic(X)$.  Choose isomorphisms
$$M_j\simeq L^{(a_{1j},\ldots,a_{rj})}, 
\quad i=1,\ldots,r, A=(a_{ij}) $$
which naturally induce isomorphisms
$$
\Gamma(M^{\nu})\simeq \Gamma(L^{A\nu}), 
\quad A\nu=(a_{11}\nu_1+\ldots+a_{1r}\nu_r,
\ldots,a_{r1}\nu_1+\ldots+a_{rr}\nu_r).
$$
Thus we find $R(X,L_1,\ldots,L_r)\simeq R(X,M_1,\ldots,M_r)$.
\end{proof}

As $\Cox(X)$ is graded by $\Pic(X)$, a free abelian group 
of rank $r$, the torus
$$T(X):=\Hom(\Pic(X),\bG_m)$$
acts on $\Cox(X)$.  Indeed, each $\nu\in \Pic(X)$ naturally yields
a character $\chi_{\nu}$ of $T(X)$, and the action is given by
$$
t\cdot \xi=\chi_{\nu}(t)\xi, \quad 
\xi \in \Cox(X)_{\nu}, t\in T(X).
$$
Thus the isomorphism constructed in Proposition~\ref{prop:unique}
is {\em not} canonical:  Two such isomorphisms differ by
the action of an element of $T(X)$.  It is precisely
this ambiguity that makes descending the
universal torsor to nonclosed fields an interesting question.

\

The following conjecture is a special case of
2.14 of \cite{HK}:

\begin{conj}[Finiteness of Cox ring]
Let $X$ be a log Fano variety. 
Then $\Cox(X)$ is finitely generated. 
\end{conj}

\begin{rem}
\label{rem:cox-fg}
Note that if $\Cox(X)$ is finitely generated 
it follows trivially that $\NE_{n-1}(X)$ is finitely generated. 
Moreover, the nef cone $\NM^1(X)$ is also finitely generated. 

Indeed, the nef cone corresponds to one of the chambers in the 
group of characters of $T(X)$ governed by the stability
conditions for points  $v\in \Spec(\Cox(X))$.  These
chambers are bounded by finitely many hyperplanes (see Theorem 0.2.3 in 
\cite{DH} for more details).  
\end{rem}

It has been conjectured by Batyrev \cite{Ba}
that the pseudo-effective cone of a Fano variety
is finitely generated. 
However, the finiteness
of the Cox ring is not a formal consequence of the
finiteness of the pseudo-effective cone.

\begin{exam}
\label{exam:9}
Let $p_1,\ldots,p_9\in H\subset \bP^3$ be nine
distinct coplanar points given as a
complete intersection of two generic
cubic curves in the hyperplane $H$, and let $X$ be the blow-up of
$\bP^3$ at these points.  Then $\NE^1(X)$ is finitely generated
but $\Cox(X)$ is not.  Indeed, $X$ is an equivariant
compactification of the additive group $\bG_a^3$,
acting by translation on the
affine space $\bP^3-H$.  The group action can
be used to show that $\NE^1(X)$ is generated
by the boundary components (see \cite{HT1}). Similarly,
one can show that the cone
$\NE_1(X)$ is generated by classes of curves in the boundary
components, e.g., the proper transform $\tilde{H}\subset X$
of $H$.  It is well-known that $\NE_1(\tilde{H})$ is infinite
\cite{KM} \S 1.23(4):  The pencil of cubic plane curves
with base locus $p_1,\ldots, p_9$ induces an elliptic fibration,
$$
\tilde{H} \ra \bP^1,
$$
for which the nine exceptional curves of $\tilde{H}\ra H$ are sections.
Addition in the group law gives an infinite number of
sections, which are also $(-1)$-curves and generators of
$\NE_1(\tilde{H})$.  These are also generators of
$\NE_1(X)$, since the sections (other than the nine
exceptional curves) intersect $\tilde{H}$ negatively.
It follows that $\NE_1(X)$ and $\NM^1(X)$ are not finitely generated
and hence $\Cox(X)$ is not finitely generated
(see Remark~\ref{rem:cox-fg}).
\end{exam}

\

\begin{prop}
\label{prop:semiample}
Let $X$ be a nonsingular projective variety
whose anticanonical divisor $-K_X$ is
nef and big.  Suppose that $D$ is a nef divisor on $X$.
Then $H^i(X,\cO_X(D))=0$ for each $i>0$ and $D$
is semiample.
\end{prop} 

\begin{proof}
The first assertion is a consequence of
Kawamata-Viehweg vanishing \cite{KM} \S 2.5.
The second is a special case of the Kawamata
Basepoint-freeness Theorem \cite{KM} \S 3.2.
\end{proof}

Proposition~\ref{prop:semiample}
largely determines the Hilbert function of the Cox ring:

\begin{coro}
\label{coro:hilbert}
Retain the assumptions of Proposition~\ref{prop:semiample}.
Then for nef classes $\nu$ we have
$$
\dim \Cox(X)_{\nu}=\chi(\cO_X(\nu)).
$$
\end{coro}

\begin{rem}
In practice, this will help us to find generators of $\Cox(X)$. 
\end{rem}

\section{Generalities on toric varieties}
\label{sect:toric}

We recall quotient constructions
of toric varieties, following Brion-Procesi \cite{brionproc},
Cox \cite{Cox}, and Thaddeus \cite{thad}.

Let $T\simeq \bG_m^r$ be a torus with character
group $\cX^*(T)$.  Suppose that $T$ acts
faithfully on the polynomial 
ring $k[x_1,\ldots,x_{n+r}]$ by the formula
$$t(x_j)=\chi_j(t)x_j, \quad t\in T, $$
where $\{\chi_1,\ldots,\chi_{n+r} \}\subset \cX^*(T).$ 
Define $M$ as the kernel of the surjective morphism
$$\chi:=(\chi_1,\ldots,\chi_{n+r}): \bZ^{n+r} \ra
\cX^*(T).$$
We interpret $M$ as the character group of the quotient torus
$\bG_m^{n+r}/T$. 
Set $N=\Hom(M,\bZ)$ so that dualizing gives
$$(\bZ^{n+r})^* \ra N \ra 0.$$
Let $e_1,\ldots,e_{n+r}$ and $e_1^*,\ldots,e_{n+r}^*$ denote
the coordinate vectors in $\bZ^{n+r}$ and  $(\bZ^{n+r})^*$;
let $\bar{e}_1^*,\ldots,\bar{e}_{n+r}^*\in N$ denote 
the images of the $e_i^*$ in $N$.  Concretely, the $\bar{e}_i^*$
are the columns of the $n\times (n+r)$ matrix of
dependence relations among the $\chi_j$.  

Consider a toric $n$-fold $X$ associated with a fan 
having one-skeleton 
$\{\bar{e}_1^*,\ldots,\bar{e}_{n+r}^*\}$.  
In particular, we assume that none of $\bar{e}_i^*$ is
zero or a positive multiple of any of the others.  
The variety $X$ is a categorical quotient of an invariant
open subset $U\subset \bA^{n+r}$ under the action of
$T$ described above (see \cite{Cox} 2.1).  
Elements $\nu\in \cX^*(T)$ classify
$T$-linearied invertible sheaves $\cL_{\nu}$ on $\bA^{n+r}$ and 
$$\Gamma(\bA^{n+r},\cL_{\nu})\simeq k[x_1,\ldots,x_{n+r}]_{\nu}.$$
We have $\rA_{n-1}(X)\simeq \cX^*(T)$ and we can identify
$$
\Gamma(\cO_X(D))\simeq k[x_1,\ldots,x_n]_{\nu(D)},
$$
where $\nu(D) \in \cX^*(T)$ is associated with
the divisor class of $D$.  The variables $x_i$ are
associated with the irreducible torus-invariant 
divisors $D_i$ on $X$ (see \cite{fulton} \S 3.4), and
the cone of effective divisors $\NE_{n-1}(X)$ is generated by 
$\{D_1,\ldots,D_{n+r}\}$.
Geometrically, the effective cone in $\cX^*(T)$ is the
image of the standard simplicial cone generated
by $e_1,\ldots,e_{n+r}$ under the projection homomorphism
$\chi:\bZ^{n+r} \ra \cX(T)$.

Recall that the {\em moving cone}
$$
\Move(X) \subset \NE_{n-1}(X)
$$
is defined as the smallest
closed subcone containing the effective divisors on $X$
without fixed components.  
\begin{prop}\label{prop:toricmove}
Retaining the notation and assumptions above,
$$
\Move(X) = \bigcap_{i=1,\ldots,n+r} 
\Cone(\chi_1,\ldots,\chi_{i-1},\chi_{i+1},\ldots,\chi_{n+r})
$$
and has nonempty interior.
\end{prop}
\begin{proof}
The fixed components of
$\Gamma(X,\cO_X(D))$ are necessarily invariant under the
torus action, hence are taken from 
$\{D_1,\ldots,D_{n+r}\}$.  Moreover, $D_i$ is fixed
in each $\Gamma(X,\cO_X(dD)),d>0$ if and only if 
$x_i$ divides each element of $k[x_1,\ldots,x_{n+r}]_{d\nu(D)}$.  
This is the case exactly when
$$\nu(D) \in \Cone(\chi_1,\ldots,\chi_{n+r}) - 
\Cone(\chi_1,\ldots,\chi_{i-1},\chi_{i+1},\ldots,\chi_{n+r}).$$

Suppose that the interior of the moving cone is empty.
After permuting indices there are two possibilities:  Either 
$\Cone(\chi_2,\ldots,\chi_{n+r})$ has
no interior, or the cones 
$\Cone(\chi_2,\ldots,\chi_{n+r})$ and
$\Cone(\chi_1,\chi_3,\ldots,\chi_{n+r})$ have nonempty interiors
but meet in a cone with positive codimension.  As the $T$-action is faithful,
the $\chi_i$ span $\cX^*(T)$.  In the first case,
$\chi_2,\ldots,\chi_{n+r}$ span a codimension-one
subspace of $\cX^*(T)$ that does not contain $\chi_1$, so 
that each dependence relation 
$$
c_1\chi_1 + \ldots + c_{n+r}\chi_{n+r}=0
$$
has $c_1=0$.  This translates into $\bar{e}_1^*=0$,
a contradiction.  In the second case, $\chi_3,\ldots,\chi_{n+r}$
span a hyperplane, and $\chi_1$ and $\chi_2$ are on opposite
sides of this hyperplane.  Putting the dependence relations 
among the $\chi_i$ in row echelon form, we obtain a unique
relation with nonzero first and second entries, and
these two entries are both positive.  This translates into
the proportionality of $\bar{e}_1^*$ and $\bar{e}_2^*$.
\end{proof}

We now seek to characterize the projective toric
$n$-folds $X$ with one-skeleton 
$\{\bar{e}_1^*,\ldots,\bar{e}_{n+r}^*\}$.
These are realized as Geometric Invariant Theory
quotients $\bA^{n+r}//T$ associated with the various 
linearizations of our $T$-action.  
We consider the graded ring
$$
R:=\sum_{d\ge 0}\Gamma(\bA^{n+r},\cL_{d\nu})=
\sum_{d\ge 0 }k[x_1,\ldots,x_{n+r}]_{d\nu}.
$$

\begin{prop}[see \cite{thad} \S 2,3]
\label{prop:whenprojective}
Retain the notation above
and set $X:=\Proj(R)$.

\begin{enumerate}
\item{
$X$ is projective over $k$ if and only if
$0$ is not contained in the convex hull of
$\{\chi_1,\ldots,\chi_{n+r} \}$.}
\item{
$X$ is toric of dimension $n$ if
$\nu$ is in the interior of 
$$
\Cone(\chi_1,\ldots,\chi_{n+r}).
$$}
\item{In this case, 
the one-skeleton of $X$ is contained in
$\{\bar{e}_1^*,\ldots,\bar{e}_{n+r}^*\}$.  
Equality holds if $\nu$ is in the interior of the moving cone 
$$\bigcap_{i=1,\ldots,n+r}
\Cone(\chi_1,\ldots,\chi_{i-1},\chi_{i+1},\ldots,\chi_{n+r}).$$}
\end{enumerate}
\end{prop}
\begin{rem}
Our proof will show that
$X$ may still be of dimension $n$ even when $\nu$
is contained in a facet of 
$$
\Cone(\chi_1,\ldots,\chi_{n+r}).
$$

Similary, the one-skeleton of $X$ may still be 
$\{\bar{e}_1^*,\ldots,\bar{e}_{n+r}^*\}$
even when $\nu$ is contained in a facet of
$$\bigcap_{i=1,\ldots,n+r}
\Cone(\chi_1,\ldots,\chi_{i-1},\chi_{i+1},\ldots,\chi_{n+r}).$$
\end{rem}
\begin{proof}
The monomials which appear in $R$ are in one-to-one correspondence
to solutions of
$$a_1\chi_1+\ldots+a_{n+r}\chi_{n+r}=d\nu, \quad
a_i \in \bZ_{\ge 0}.$$
In geometric terms, the monomials appearing
in $R$ coincide with the elements of $\bZ^{n+r}$
in the cone
$$\chi^{-1}(\Cone(\nu))\cap \Cone(e_1,\ldots,e_{n+r}).$$
By Gordan's Lemma in convex geometry, $R$ is generated
as a $k$-algebra by a finite set of monomials 
$x^{m_1},\ldots,x^{m_s}$.
The monomials appearing in the $d$th graded piece $R_d$ coincide
with elements of $\bZ^{n+r}$ in the polytope
$$P_{d\nu}:=\chi_{\bR}^{-1}(d\nu) \cap \Cone(e_1,\ldots,e_{n+r}).$$
Note that $\chi^{-1}(d\nu)$ is a translate of $M$.

For the first part, recall that 
$\Proj(R)$ is projective over $\Spec(R_0)$,
where $R_0$ is the degree-zero part.  
Now $0$ is in the convex hull of $\{\chi_1,\ldots,\chi_{n+r}\}$
if and only if there are nonconstant elements of $R$ of degree
zero.  Our hypothesis just says that $R_0=k$ and thus is equivalent
to the projectivity of $X$ over $k$.  

As for the second part,  $T$ acts 
on $R$ by homotheties and thus acts trivially on $\Proj(R)$,
so we have an induced action of $\bG_m^{n+r}/T$ on $\Proj(R)$.
We claim this action is faithful, so 
the quotient is toric of dimension $n$.  
Let $\mu_1,\ldots,\mu_n$ be generators for 
$M=\cX^*(\bG_m^{n+r}/T)$. 
Choose $v \in \bZ^{n+r}$ in the interior of 
$\Cone(e_1,\ldots,e_{n+r})$ so that $\chi_{\bR}(\Cone(v))=\Cone(\nu)$.
Replacing $v$ by a suitably large integral multiple, 
we may assume each $v+\mu_i,i=1,\ldots,n$, is in 
$\Cone(e_1,\ldots,e_{n+r})$.  If $\chi(v)=d\nu$ then $R_d$ contains
a set of generators for $M$, so the induced representation
of $\bG_m^{n+r}/T$ on $R_d$ is faithful.  

For the third part, we extract the fan classifying $X$
from $P_{d\nu}$, following \cite{fulton} \S 1.5 and 3.4:
For each face $Q$ of $P_{d\nu}$, consider the cone
$$\sigma_Q=\{v \in N_{\bR}: \left<u,v \right> \leq \left<u',v \right>
\text{ for all }u \in Q, u' \in P_{d\nu} \}.$$
This assignment is inclusion reversing, so the one-dimensional
cones of the fan correspond to facets of $P_{d\nu}$.  Moreover,
each facet of $P_{d\nu}$ is induced by one of the facets
of $\Cone(e_1,\ldots,e_{n+r})$.  The corresponding one-dimensional
cone in $N_{\bR}$ is spanned by $\bar{e}_i^*$.  
It remains to verify that each facet of $\Cone(e_1,\ldots,e_{n+r})$
actually induces a facet of $P_{d\nu}$.  
The hypothesis that $\nu$ is in the moving cone 
means that $P_{d\nu}$ intersects each of the 
$\Cone(e_1,\ldots,e_{i-1},e_{i+1},\ldots,e_{n+r})$.
If $\nu$ is in the interior of the moving cone then
the intersection of $P_{d\nu}$ with 
$\Cone(e_1,\ldots,e_{i-1},e_{i+1},\ldots,e_{n+r})$
meets the relative interior of this cone, hence
this cone induces a facet of $P_{d\nu}$.  
\end{proof}

Proposition \ref{prop:whenprojective} yields the following
nice consequence:

\begin{prop} 
\label{prop:projmodels}
Let $X$ be a complete toric variety and $\nu$
a divisor class in the interior of $\Move(X)$. 
Then there exists a projective toric variety $Y_{\nu}$,
with the same one-skeleton as $X$, and polarized by $\nu$.  
\end{prop}

\noindent
For generic $T$-linearized invertible sheaves on $\bA^{n+r}$,
all semistable points are actually stable;  hence 
$Y_{\nu}$ is a simplicial toric
variety for generic $\nu$ (see \cite{brionproc} 1.2 and \cite{Cox} 2.1).  
For the special values $\nu_0$, contained in the
walls of the chamber decomposition of \cite{thad}, this
fails to be the case.  However, for each special $\nu_0$,
there exists a generic $\nu$ so that 
$\Cone(\nu)$ is very close to $\Cone(\nu_0)$ and
there is a projective, torus-equivariant morphism
$Y_{\nu} \ra Y_{\nu_0}$ \cite{thad} 3.11.  
The polarization associated to $\nu_0$ pulls back to $Y_{\nu}$,
so we obtain the following:

\begin{prop} 
\label{prop:goodmodels}
Let $X$ be a complete toric variety and $\nu_0$
a divisor class in the moving cone of $X$. 
Then there exists a simplicial projective toric variety
$Y$, with the same one-skeleton as $X$, so that $\nu_0$ is
semiample on $Y$.
\end{prop}
Of course, $\nu_0$ is big when it is in the interior
of the effective cone.

\section{The $\Esix$ cubic surface}
\label{sect:e6}

By definition, 
the {\em $\Esix$ cubic surface} is given by the homogeneous
equation
\begin{equation} 
\label{eqn:main}
S=\{(w,x,y,z):xy^2+yw^2+z^3=0\}\subset \bP^3.
\end{equation}
We recall some elementary properties (see \cite{bruce-wall}
for more details on singular cubic surfaces):

\begin{prop}
\label{prop:classical}
$
{}
$
\begin{enumerate}
\item The surface $S$ has a single singularity at the
point $p:=(0,1,0,0)$, of type $\Esix$.
\item $S$ is the unique cubic surface with this property,
up to projectivity.
\item $S$ contains a unique line, satisfying the equations
$y=z=0$.
\end{enumerate}
\end{prop}

Any smooth cubic surface may be
represented as the blow-up of $\bP^2$ at six points
in `general position'.  There is an analogous property
of the $\Esix$ cubic surface:

\begin{prop}
\label{prop:sixpoints}
The $\Esix$ cubic surface $S$ is 
the closure of the image of $\bP^2$ under the
linear series
$$w=a^2c \quad x=-(ac^2+b^3) \quad y=a^3 \quad z=a^2b,$$
where 
$$
\Gamma(\bP^2,\cO_{\bP^2}(1))=\langle a,b,c\rangle.
$$
This map is the inverse of the projection of $S$ 
from the double point  $p$.
The affine open subset
$$
\bA^2:=\{a\neq 0\} \subset \  \bP^2 
$$
is mapped isomorphically onto $S-\ell$.  In particular,
$S\setminus \ell \simeq \bA^2$, so the $\Esix$ cubic surface is a 
compactification of $\bA^2$.
\end{prop}

\begin{rem}
Note that $S$ is not an {\em equivariant} compactification of 
$\bG_a^2$, so the general theory of \cite{CLT} does
not apply.

Indeed, if $S$ were an equivariant compactification of $\bG_a^2$
then the projection from $p$ would be $\bG_a^2$-equivariant (see \cite{HT1}).
Therefore, the map $\bP^2 \dashrightarrow S$ given above 
has to be $\bG_a^2$-equivariant. 
The only $\bG_a^2$-action on $\bP^2$ under which a line is invariant
is the standard translation action \cite{HT1}.  
However, the linear series above is not invariant under the 
standard translation action 
$$
b\mapsto b+\beta a \quad c \mapsto c+\gamma a.
$$
\end{rem}

\

We proceed to compute the effective cone of the  minimal 
resolution $\phi_{\ell}:\tS\ra S$. 
Let $\ell \subset \tS$ be the proper transform of the
line mentioned in Proposition~\ref{prop:classical}.  

\begin{prop}
\label{prop:PIC}
The Picard group $\Pic(\tS)$ is a free abelian group of rank seven, generated
by $\ell$ and the exceptional curves of $\phi_{\ell}$.  
For a suitable ordering $\{F_1,F_2,F_3,F_4,F_5,F_6\}$
of the exceptional curves, the intersection
pairing takes the form
\begin{equation}
\label{Fmatrix}
\begin{array}{c|ccccccc}
 & F_1 & F_2 & F_3 & \ell & F_4 & F_5 & F_6 \\
\hline
F_1 & -2 & 0 & 1 & 0 & 0 & 0 & 0 \\
F_2 & 0 & -2 & 0 & 0  & 0 & 0 & 1   \\
F_3 & 1 & 0 & -2 & 0 & 0 & 0 & 1  \\
\ell& 0 & 0 & 0 & -1  & 1 & 0  & 0   \\
F_4 & 0 & 0 & 0 & 1 & -2 & 1 & 0   \\
F_5 & 0 & 0 & 0 & 0 & 1 & -2 & 1    \\
F_6 & 0 & 1 & 1 & 0 & 0 & 1 & -2  
\end{array}.
\end{equation}
\end{prop}

\begin{figure}[t]
\label{E6figure}
\centerline{\includegraphics{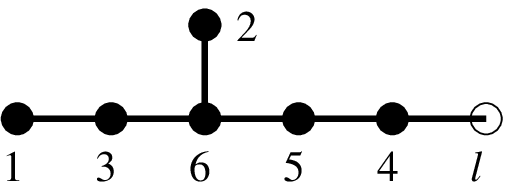}}
\caption{Dynkin diagram of $\Esix$}
\end{figure}

\begin{prop}
\label{prop:effcone}
The effective cone $\NE(\tS)$ 
is simplicial and generated by 
$\Phi:=\{F_1,F_2,F_3,\ell,F_4,F_5,F_6\}$.  
Each nef divisor is contained in the monoid generated by 
the divisors:
\begin{eqnarray*}
A_1      &=&\quad \quad  F_2 + F_3 + 2\ell + 2F_4 + 2F_5 +  2F_6  \\
A_2      &=& F_1  + F_2 + 2F_3 + 3\ell + 3F_4 + 3F_5 +  3F_6   \\
A_3      &=& F_1 + 2F_2 + 2F_3 + 4\ell + 4F_4 + 4F_5 +  4F_6  \\
A_{\ell} &=& 2F_1 + 3F_2 + 4F_3 + 3\ell + 4F_4 + 5F_5 +  6F_6  \\
A_4      &=& 2F_1 + 3F_2 + 4F_3 + 4\ell + 4F_4 + 5F_5 +  6F_6   \\
A_5      &=& 2F_1 + 3F_2 + 4F_3 + 5\ell + 5F_4 + 5F_5 +  6F_6  \\
A_6      &=& 2F_1 + 3F_2 + 4F_3 + 6\ell + 6F_4 + 6F_5 +  6F_6  
\end{eqnarray*}
Moreover $A_{\ell}$ is the anticanonical class $-K_{\tS}$
and its sections induce the resolution morphism $\phi_{\ell}:\tS \ra S$.  
\end{prop}

\begin{proof}
The intersection form in terms of
$A:=\{A_1,\ldots,A_6  \}$ is:
\begin{equation}
\label{Amatrix}
\begin{array}{c|ccccccc}
 & A_1 & A_2 & A_3 & A_{\ell} & A_4 & A_5 &  A_6  \\
\hline
A_1      & 0 & 1 & 1 & 2 & 2 & 2 &  2  \\
A_2      & 1 & 1 & 2 & 3 & 3 & 3 &  3   \\
A_3      & 1 & 2 & 2 & 4 & 4 & 4 &  4  \\
A_{\ell} & 2 & 3 & 4 & 3 & 4 & 5 &  6  \\
A_4      & 2 & 3 & 4 & 4 & 4 & 5 &  6   \\
A_5      & 2 & 3 & 4 & 5 & 5 & 5 &  6  \\
A_6      & 2 & 3 & 4 & 6 & 6 & 6 &  6  
\end{array}.
\end{equation}
This is the inverse of the intersection
matrix \eqref{Fmatrix} written in terms of the basis $\Phi$,
so the $A_i$ generate the dual to 
$\Cone(\Phi)$.  Observe that all the 
entries of matrix \eqref{Amatrix} are nonnegative and
$$\Cone(A)\subset \Cone(\Phi).$$

Suppose that $D$ is an effective divisor on $\tS$.  
We write $D$ as a sum of the fixed components
contained in $\{F_1,\ldots,F_6,\ell\}$ and the
parts moving relative to $\Phi$:
$$
D=M_{\Phi}+a_1F_1+\ldots +a_6F_6+a_{\ell}\ell, \quad
a_1,\ldots,a_6,a_{\ell}\ge 0.
$$
{\em A priori}, $M_{\Phi}$ may have fixed components,
but they are not contained in $\Phi$ (however, see Lemma~\ref{lemm:fin}).  
It follows that 
$M_{\Phi}$ intersects each element of $\Phi$ nonnegatively,
i.e., it is contained in $\Cone(A)$ and thus
in $\Cone(\Phi)$.  We conclude that $D\in \Cone(\Phi)$.
Since $A_1, ..., A_6, A_{\ell}$ generate $\rN_1(\tS)$ over $\bZ$
each nef divisor can be written as a nonnegative linear combination
of these divisors. 

To see that $A_{\ell}$ is the anticanonical divisor, we apply
adjunction
$$
K_{\tS}F_i=0, i=1,\ldots, 6 \quad K_{\tS}\ell=-1.
$$
Nondegeneracy of the intersection form implies
$A_{\ell}=-K_{\tS}$.  Since $S$ has rational double points,
the resolution map $\phi_{\ell}$ is crepant, i.e., 
$\phi_{\ell}^*K_S=K_{\tS}$.
Thus 
$$
\Gamma(A_{\ell})=\Gamma(-K_{\tS})=\Gamma(-\phi_{\ell}^*K_S)=
\Gamma(\phi_{\ell}^*\cO_S(+1))
$$
so the sections of $A_{\ell}$ induce $\phi_{\ell}$.  
\end{proof}

Choose nonzero sections
$\xi_1,\ldots,\xi_{\ell}$
generating $\Gamma(F_1),\ldots,\Gamma(\ell)$:
$$\Gamma(F_1)=\left<\xi_1\right>, \ldots, 
\Gamma(F_6)=\left<\xi_6\right>,  \Gamma(\ell)=\left<\xi_{\ell}\right>.$$
These are canonical up to scalar multiplication.
Each effective divisor
$$D=b_1F_1+b_2F_2+b_3F_3+b_{\ell}\ell + b_4F_4+ b_5F_5+b_6F_6$$
has a distinguished nonzero section
$$\xi^{(b_1,b_2,b_3,b_{\ell},b_4,b_5,b_6)}:=
\xi_1^{b_1}\ldots \xi_6^{b_6}\xi_{\ell}^{b_{\ell}}.$$ 
The distinguished section of $A_j$
is denoted $\xi^{\alpha(j)}$.  
Note that we have an injective ring homomorphism 
\begin{equation}
\label{eqn:ring-h}
k[\xi_1,\ldots, \xi_6,\xi_{\ell}]\ra \Cox(\tilde{S}). 
\end{equation}

There is a partial order on the monoid of
effective divisors of $\tS$:  $D_1\prec D_2$
if $D_2-D_1$ is effective.  The restriction of this
order to the generators of the nef cone
is illustrated in the diagram below:
$$
\begin{array}{ccccc}
 & & A_6 & & \\
 & & | & & \\
 & & A_5 & & \\
 & & | & & \\
 & & A_4 & & \\
 & \diagup &  & \diagdown & \\
A_{\ell} &  &  &  & A_3 \\
 & \diagdown &  & \diagup & \\
 & & A_2 & & \\
 & & | & & \\
 & & A_1 & & 
\end{array}
$$
Whenever $D_1\prec D_2$
we have an inclusion 
$$\Gamma(D_1)\hra \Gamma(D_2)$$
which is natural up to scalar multiplication.
Indeed, express
$$D_1-D_2=b_1F_1+b_2F_2+\ldots+b_6F_6+b_{\ell}\ell,
\quad b_j \ge 0$$
so we have
\begin{eqnarray*}
s_1 &\mapsto&
\xi^{(b_1,b_2,b_3,b_{\ell},b_4,b_5,b_6)}s_1 \\
\Gamma(D_1) & \hra & \Gamma(D_2).
\end{eqnarray*}

The homomorphism \eqref{eqn:ring-h} is not surjective, and we
now look for generators of 
$\Cox(\tS)$ beyond the $\xi_j$.  
Consider the subring
$$
\Cox_a(\tS)=\bigoplus_{\nu\in \NM(\tS)}\Cox(\tS)_{\nu}
$$ 
obtained by restricting to degrees corresponding 
to nef classes on $\tS$.  
The following lemma implies that 
any homogeneous element $\mathsf s_D\in \Cox(\tS)$ 
can be written in the form
$$
\mathsf s_D=m_D\xi_1^{b_1}\cdots \xi_6^{b_6}\xi_{\ell}^{b_\ell}
$$
with nonnegative exponents and $m_D\in \Cox_a(\tS)$.  
\begin{lem}
\label{lemm:fin}
Let $D$ be an effective divisor on $\tS$ with fixed part $F_D$ and moving
part $M_D$. Then $F_D$ is supported in 
$\{ F_1, \ldots , F_6, \ell\}$, and $M_D$ is     
a linear combination of 
$A_1, \ldots, A_6,A_{\ell}$ with nonnegative
coefficients.    
\end{lem}
\begin{proof} 
Clearly $M_D$ is nef, so the description of the nef divisors in
Proposition~\ref{prop:effcone} gives the expression
in terms of the $A_i$.  
Proposition~\ref{prop:semiample} shows $M_D$ is
semiample with vanishing higher cohomology;  the last part of
Proposition~\ref{prop:effcone} gives the requisite
positivity of the anticanonical class.  

Let $F$ be a fixed component of $D$ not supported in 
$\{ F_1, \ldots , F_6, \ell\}$. 
To arrive at a contradiction, we need to show that
$h^0(\cO_{\tS}(M_D+F))> h^0(\cO_{\tS}(M_D))$. 
Since $M_D$ has vanishing higher cohomology and 
$$
h^2(\cO_{\tS}(F+M_D))=h^0(\cO_{\tS}(K_{\tS}-F-M_D)) =0 
$$
it suffices to show that 
$$
\chi(\cO_{\tS}(F+M_D))> \chi(\cO_{\tS}(M_D)).
$$
By Riemann-Roch, it suffices to show that 
$$
F^2+2M_DF-K_{\tS}F >0
$$
or, equivalently, 
$$
F^2+K_{\tS}F +2M_DF -2K_{\tS}F = 2g(F)-2 + 2M_DF -2K_{\tS}F >0.
$$
Since $F$ is irreducible, $g(F)\ge 0$ 
and $M_DF\ge 0$ and $-K_{\tS}F > 1$, 
as  $M_D$ is nef and $-K_{\tS}$ is nonpositive only along the
exceptional curves and has degree 1 only on the line $\ell$ 
(see Proposition~\ref{prop:classical}). 
\end{proof}

Corollary \ref{coro:hilbert} gives
the dimensions of the graded pieces of $\Cox_a(\tS)$.
We focus first on the generators of the nef cone,
introducing sections $\tau_j\in \Gamma(A_j)$ as needed to 
achieve the prescribed dimensions:
\begin{eqnarray*}
\Gamma(A_1)&=&\left<\xi^{\alpha(1)},\tau_1\right>\\
\Gamma(A_2)&=&\left<\xi^{\alpha(2)},\xi^{\alpha(2)-\alpha(1)}\tau_1,
\tau_2 \right> \\
\Gamma(A_{\ell})&=&\left<\xi^{\alpha(\ell)},\xi^{\alpha(\ell)-\alpha(1)}\tau_1,
\xi^{\alpha(\ell)-\alpha(2)}\tau_2,\tau_{\ell} \right> 
\end{eqnarray*}
The sections of $A_{\ell}$ induce $\phi_{\ell}:\tS\ra S\subset \bP^3$
by Proposition \ref{prop:effcone}, and can be identified with
the coordinates $w,x,y,z$ of Equation~\eqref{eqn:main}.  
Since $A_1\prec A_2 \prec A_{\ell}$,
we have 
$$\Gamma(A_1)\hookrightarrow \Gamma(A_2) \hookrightarrow\Gamma(A_{\ell}).$$
We can identify $\Gamma(A_1)=\left<y,z\right>$;  these
correspond to projecting $S$ from the line $\ell=\{y=z=0\}$
and induce a conic bundle structure
$$
\phi_1:\tS \ra \bP^1.
$$
We have $\Gamma(A_2)=\left<x,y,z\right>$;  these correspond
to projecting $S$ from the singularity $p=\{w=y=z=0\}$ and induce the
blow-up realization
$$
\phi_2:\tS \ra \bP^2.
$$
Therefore, we may choose $\tau_1,\tau_2,$ and $\tau_{\ell}$ so that 
$$
y=\xi^{\alpha(\ell)} \quad w=\xi^{\alpha(\ell)-\alpha(2)}\tau_2 \quad
z=\xi^{\alpha(\ell)-\alpha(1)}\tau_1 \quad x=\tau_{\ell}.
$$

We obtain the following induced sections for $A_3,A_4,A_5,$ and $A_6$:
\begin{eqnarray*}
\Gamma(A_3)&=&\left<\xi^{\alpha(3)},\xi^{\alpha(3)-\alpha(1)}\tau_1,
\xi^{\alpha(3)-\alpha(2)}\tau_2,\xi^{\alpha(3)-2\alpha(1)}\tau_1^2\right>\\
\Gamma(A_4)&=&\left<\xi^{\alpha(4)},\xi^{\alpha(4)-\alpha(1)}\tau_1,
\xi^{\alpha(4)-\alpha(2)}\tau_2,
\xi^{\alpha(4)-\alpha(\ell)}\tau_{\ell}, 
\xi^{\alpha(4)-2\alpha(1)}\tau_1^2\right>\\
\Gamma(A_5)&=&\left<\xi^{\alpha(5)},\xi^{\alpha(5)-\alpha(1)}\tau_1,
\xi^{\alpha(5)-\alpha(2)}\tau_2,
\xi^{\alpha(5)-\alpha(\ell)}\tau_{\ell}, 
\xi^{\alpha(5)-2\alpha(1)}\tau_1^2, \right. \\
& &\quad \quad \left. \xi^{\alpha(5)-\alpha(1)-\alpha(2)}\tau_1\tau_2 \right>\\
\Gamma(A_6)&=&\left<\xi^{\alpha(6)},\xi^{\alpha(6)-\alpha(1)}\tau_1,
\xi^{\alpha(6)-\alpha(2)}\tau_2,
\xi^{\alpha(6)-\alpha(\ell)}\tau_{\ell}, 
\xi^{\alpha(6)-2\alpha(1)}\tau_1^2, \right. \\
& & \quad \quad \left.
\xi^{\alpha(6)-\alpha(1)-\alpha(2)}\tau_1\tau_2,
\xi^{\alpha(6)-2\alpha(2)}\tau^2_2,
\xi^{\alpha(6)-3\alpha(1)}\tau_1^3\right>
\end{eqnarray*}
Equation~\eqref{eqn:main} gives the relation
$$
\tau_{\ell}\xi^{2\alpha(\ell)}+\tau^2_2\xi^{3\alpha(\ell)-2\alpha(2)}+
\tau_1^3\xi^{3\alpha(\ell)-3\alpha(1)}=0.
$$
Dividing by a suitable monomial $\xi^{\beta}$, we obtain
$$
\tau_{\ell}\xi_{\ell}^3\xi_4^2\xi_5 + \tau_2^2  \xi_2 + 
\tau_1^3 \xi_1^2 \xi_3=0,
$$
a dependence relation in $\Gamma(A_6)$.
This is the only such relation:
Any other relation,
after multiplying through by $\xi^{\beta}$, yields a cubic 
form vanishing on $S\subset \bP^3$,
but equation ~\eqref{eqn:main} is the only such form.
It follows that the sections given above for $A_1,\ldots,A_5$
form bases for $\Gamma(A_1),\ldots,\Gamma(A_5)$.  

Since 
$$A_3 \prec A_4 \prec A_5 \prec A_6 \prec 2A_{\ell}$$
we have
\begin{eqnarray*}
\Gamma(A_3)&\hookrightarrow &\Gamma(A_4) \hookrightarrow
\Gamma(A_5) \hookrightarrow \Gamma(A_6) \\
& \hookrightarrow &\Gamma(2A_{\ell})=
\left<w^2,wx,wy,x^2,xy,xz,y^2,yz,z^2 \right>
\end{eqnarray*}
and identifications
\begin{eqnarray*}
\Gamma(A_3)&=&\left<y^2,yz,
wy,z^2 \right>\\
\Gamma(A_4)&=&\left<y^2,yz,
wy,xy,z^2\right>\\
\Gamma(A_5)&=&\left<y^2,yz,wy,xy,z^2,wz\right> \\
\Gamma(A_6)&=&\left<y^2,yz,wy,xy,z^2,wz,w^2\right>. 
\end{eqnarray*}
The sections of $A_3$ induce a morphism
$$\phi_3:\tilde{S} \ra  \bP^3$$
onto a quadric surface with a single ordinary double
point.  The sections of $A_4$ induce
a morphism 
$$\phi_4:\tilde{S} \ra \bP^4$$
with image a quartic Del Pezzo surface with
a rational double point of type $\mathbf{D}_5$.  The sections of
$A_5$ induce a morphism
$$\phi_5:\tilde{S} \ra \bP^5$$ 
with image a quintic Del Pezzo surface with a rational double point
of type $\mathbf{A}_4$.  The sections of $A_6$
induce a morphism 
$$\phi_6:\tilde{S} \ra \bP^6$$
with image a sextic Del Pezzo surface with two rational double points,
of types $\mathbf{A}_1$ and $\mathbf{A}_2$.  

We summarize this analysis in the following proposition
\begin{prop}\label{prop:semiampsect}
Every section of $A_j, j=1,2,3,\ell,4,5,6,$
can be expressed as a polynomial in 
$\xi_1,\ldots,\xi_6,\xi_{\ell},\tau_1,\tau_2,\tau_6$.  The only
dependence relation among these is
$$
\tau_{\ell}\xi_{\ell}^3\xi_4^2\xi_5 + \tau_2^2  \xi_2 + 
\tau_1^3 \xi_1^2 \xi_3=0
$$
in $\Gamma(A_6)$.  Each $A_j$ is globally generated and induces
a morphism 
$$\phi_j:\tS \ra \bP^{\chi-1}, \quad \chi=\chi(\cO_{\tS}(A_j)).$$
\end{prop}

\

The remainder of this section is devoted to proving the following:
\begin{theo}
\label{theo:main}
The homomorphism 
$$
\varrho\,:\, 
k[\xi_1,...,\xi_6,\xi_{\ell}, \tau_1, \tau_2,\tau_\ell] /
\langle \tau_{\ell}\xi_{\ell}^3\xi_4^2\xi_5 + \tau_2^2  \xi_2 + 
\tau_1^3 \xi_1^2 \xi_3\rangle \ra \Cox(\tS)
$$
is an isomorphism. 
\end{theo}

If $\varrho$ were not injective, its kernel would have
nontrivial elements in degree $\nu=dA_{\ell}$, for some $d$
sufficiently large.  These translate
into homogeneous polynomials of degree $d$ vanishing on
$S\subset \bP^3$.  All such polynomials are multiples of the
cubic form defining $S$, which itself is a multiple of the relation
we already have. 

It remains to show that $\varrho$ is surjective.  By
Proposition~\ref{prop:effcone}, Lemma~\ref{lemm:fin} 
and the analysis of the sections of the $A_i$,
it suffices to prove:
\begin{prop}\label{prop:generate}
$\varrho$ is surjective in degrees corresponding to nef
divisor classes of $\tS$.
\end{prop}
\begin{lem}\label{lem:embed}
For any positive integers 
$c_1,c_2,c_3,c_{\ell},c_4,c_5,c_6$, the image of 
$$\Gamma(A_1)^{\otimes c_1} \otimes \ldots \otimes
  \Gamma(A_{\ell})^{\otimes c_{\ell}} \otimes
\ldots  \otimes  \Gamma(A_6)^{\otimes c_6}
\longrightarrow \Gamma(c_1A_1+\ldots+c_6A_6)$$
is a linear series embedding $\tS$.
\end{lem}
\begin{proof}
Proposition~\ref{prop:semiampsect} says that each $A_j$ is globally
generated, so if the image of 
$$\Gamma(A_1) \otimes \ldots \otimes \Gamma(A_6)
\longrightarrow \Gamma(A_1+\ldots+A_6)$$
embed $\tS$ then the general result follows.
We use the standard criterion:
a linear series gives an embedding iff any length-two subscheme 
$\Sigma \subset \tS$ imposes two independent conditions on the
linear series.  

First, suppose the support of $\Sigma$ is not contained in the exceptional
locus of $\phi_{\ell}:\tS \ra S$, i.e., the curves $F_1,F_2,F_3,F_4,F_5,F_6$.
Then $\phi_{\ell}$ maps $\Sigma$ to a subscheme of length two, which
imposes independent conditions on $\Gamma(A_{\ell})$,
and thus independent conditions on the linear series in question.  
Second, suppose that $\Sigma \subset F_j$ for some $j$ (resp.
$\Sigma \subset \ell$).  Since $A_j \cdot F_j=1$ (resp. 
$A_{\ell}\cdot \ell=1$), $\phi_j$ maps $F_j$ (resp. $\ell$) 
isomorphically onto a line.  It follows that $\Sigma$
imposes independent conditions on $\Gamma(A_j)$.  
Third, suppose that $\Sigma$ is reduced with support in $F_i$ and $F_j$,
but is not contained in either $F_i$ or $F_j$.  Consider the
chain of rational curves containing $F_i$ and $F_j$ (see 
Figure~\ref{E6figure}.)  There exists a curve $F_k$
in this chain so that $\phi_k(F_i)\neq \phi_k(F_j)$, so
$\Sigma$ imposes independent conditions on $\Gamma(A_k)$.  
Fourth, suppose that $\Sigma$ is nonreduced and supported
in $F_j$ but not contained in any $F_i$ or $\ell$.  The
morphism $\phi_j$ ramifies at points where $F_j$ meets
one of the other exceptional curves, and the kernel of the
tangent morphism $d\phi_j$ consists of the tangent vectors
to the curves contracted by $\phi_j$.   It follows
that $\phi_j(\Sigma)$ has length two and imposes independent
conditions on $\Gamma(A_j)$.  
\end{proof}

The polynomial ring 
$$k[\xi_1,\ldots,\xi_6,\xi_{\ell},\tau_1,\tau_2,\tau_{\ell}]$$
is graded by the N\'eron-Severi group of $\tS$
$$\deg(\xi_j)=F_j, j=1,\ldots,6 \ \deg(\xi_{\ell})=\ell \ 
\deg(\tau_j)=A_j, j=1,2,\ell.$$
This gives an action of the N\'eron-Severi torus
$T(\tS)$ on the corresponding affine space $\bA^{10}$.  

We consider the projective toric varieties that arise as 
quotients of $\bA^{10}$ by $T(\tS)$.  As sketched in
\S \ref{sect:toric}, 
these varieties have the one-skeleton
$$ x_1=(0,1,2),x_2=(1,1,3),x_3=(1,2,4), x_{\ell}=(2,3,3),x_4=(2,3,4) $$
$$ x_5=(2,3,5), x_6=(2,3,6), t_1=(-1,0,0), t_2=(0,-1,0), 
t_{\ell}=(0,0,-1) $$
where the $x_j$ correspond to the $\xi_j$ and the $t_j$ correspond
to the $\tau_j$.  

\begin{lem} \label{lem:computemoving}
Let $X$ be a toric threefold with one-skeleton $\{x_1,\ldots,t_{\ell}\}$
and divisor class-group $\cX^*(T(\tS))=\rN_1(\tS)$.  
Then 
$$\Move(X)=\Cone(A_1,\ldots,A_6,A_{\ell}).$$
\end{lem}
\begin{proof}
Proposition~\ref{prop:toricmove} reduces this to computing the
intersection of the cones generated by subsets of
$$\{F_1,\ldots,F_6,\ell,A_1,A_2,A_{\ell}\}$$
with nine elements.  Since $A_1,A_2,A_{\ell}$ are effective
combinations of the classes $F_1,\ldots,F_6,$ and $\ell$, it suffices to
compute
$$
\Cone(F_1,\ldots,\hat{\ell},\ldots,A_{\ell}) \bigcap
\left(\bigcap_{i=1,\ldots,6} \Cone(F_1,\ldots,\hat{F_i},\ldots,
A_{\ell})  \right).
$$

This intersection obviously contains $A_1,A_2,$ and $A_{\ell}$,
and it is a straightforward computation to show that it also 
contains $A_3,A_4,A_5,A_6$.  For the reverse inclusion, 
suppose that $D$ is contained in the intersection.
Considering $D$ as a divisor
on $\tS$, we see that 
$$D\cdot F_1,\ldots,D\cdot F_6,D\cdot \ell$$
are all nonnegative.  Thus $D$ is an effective sum of $A_j$
by Proposition \ref{prop:effcone}.  
\end{proof}

Combining Lemmas~\ref{lem:computemoving} and \ref{lem:embed}
with Propositions~\ref{prop:whenprojective} and 
\ref{prop:semiampsect}, we obtain the following
\begin{prop}\label{prop:toricembedding}
Let $\nu$ be an ample divisor on $\tS$.  Then there exists a
projective toric variety $Y_{\nu}$ with one-skeleton
$\{x_1,\ldots,t_{\ell}\}$ and polarization $\nu$, and an embedding
$\tS \hookrightarrow Y_{\nu}$ with the following properties:
\begin{enumerate}
\item{the divisor class group of $Y_{\nu}$
is isomorphic to the divisor class group of $\tS$
so that the moving cone of $Y_{\nu}$
is identified with the nef cone of $\tS$;}
\item{
the equation for $\tS$ in the Cox ring of $Y_{\nu}$ is 
$$
\tau_{\ell}\xi_{\ell}^3\xi_4^2\xi_5 + \tau_2^2  \xi_2 +
\tau_1^3 \xi_1^2 \xi_3=0
$$
and $[\tS]=A_6$ in the divisor class group of
$Y_{\nu}$;}
\item{$\Cox(Y_{\nu})=k[\xi_1,\ldots,\tau_{\ell}]$
and is mapped isomorphically to the image of the homomorphism $\varrho$.}
\end{enumerate}
\end{prop}

For each toric variety $Y_{\nu}$, we can consider the
exact sequence of sheaves
$$0 \ra I_{\tS} \ra \cO_{Y_{\nu}}
\ra \cO_{\tS} \ra 0,$$
where $I_{\tS}\simeq \cO_{Y_{\nu}}(-A_6)$ is the ideal sheaf
of $\tS$.  Given an element $\theta$ in the divisor class group 
of $Y_{\nu}$, we can twist to obtain
$$0 \ra I_{\tS}(\theta) \ra \cO_{Y_{\nu}}(\theta)
\ra \cO_{\tS}(\theta) \ra 0.$$
We should make precise what we mean by the twist
$\cF(\theta)$ of a coherent sheaf
$\cF$ on $Y_{\nu}$:  Realize $\cF$
as the sheafification of a graded module $F$ over $\Cox(Y_{\nu})$
(which exists by \cite{Mus} Theorem 1.1, \cite{Cox} Proposition 3.1), 
shift $F$ by $\theta$,
and then resheafify the shifted module to obtain $\cF(\theta)$.  
Twisting respects exact sequences \cite{Cox} 3.1.

The anticanonical divisor of a toric variety is the 
sum of the invariant divisors \cite{fulton} p. 89, so 
$$-K_{Y_{\nu}}=F_1+\ldots +F_6+\ell+A_1+A_2+A_{\ell}=A_{\ell}+A_6$$
and we can rewrite our exact sequence as
$$0 \ra \cO_{Y_{\nu}}(K_{Y_{\nu}}+A_{\ell}+\theta) \ra \cO_{Y_{\nu}}(\theta)
\ra \cO_{\tS}(\theta) \ra 0.$$

Suppose that $\theta$ corresponds to a nef class on $\tS$;  we shall
prove that $\varrho$ is surjective in degree $\theta$,
thus proving Proposition~\ref{prop:generate} and
Theorem~\ref{theo:main}.  Since
$$\Gamma(\cO_{Y_{\nu}}(\theta))\simeq 
k[\xi_1,\ldots,\xi_6,\xi_{\ell},\tau_1,\tau_2,\tau_{\ell}]_{\theta}$$
it suffices to show that 
$$H^1(\cO_{Y_{\nu}}(K_{Y_{\nu}}+A_{\ell}+\theta))=0.$$
We apply Proposition~\ref{prop:goodmodels}, with $\nu_0=A_{\ell}+\theta$,
to get a {\em simplicial} toric variety $Y_{\nu}$ on which
$A_{\ell}+\theta$ is nef.  
As $A_{\ell}$ is in the interior
of the effective cone of $Y_{\nu}$, $A_{\ell}+\theta$ is also big.  
Note that $Y_{\nu}$ has finite-quotient singularities, which are
log terminal \cite{KM} \S 5.2.  The desired vanishing follows
from Theorem 2.17 of \cite{Ko}.  Alternately, we could apply
Theorem 0.1 of \cite{Mus}, which applies in arbitrary characteristic
and obviates the need to pass to a simplicial model.

\

\section{$\Dfour$ cubic surface} 
\label{sect:d4}

The strategy of the previous section can applied to other
surfaces as well.  Here we illustrate it in the case
of a cubic surface given by the homogeneous equation
$$
S=\{(x_1,x_2,x_3,w)\, :\,  w(x_1+x_2+x_3)^2 = x_1x_2x_3 \} \subset \bP^3.
$$
We summarize its properties: 
\begin{enumerate}
\item $S$ has a single singularity at the point $p=(0,0,0,1)$ of type $\Dfour$.
\item $S$ contains 6 lines with the equations
$$
\begin{array}{cc}
\ell_1':= \{ w=x_1=0\}  &  m_1':= \{ x_1=x_2+x_3=0\} \\ 
\ell_2':= \{ w=x_2=0\}  &  m_2':= \{ x_2=x_1+x_3=0\}  \\
\ell_3':= \{ w=x_3=0\}  &  m_3':= \{x_3=x_1+x_2=0\} 
\end{array}
$$
\item $S$ is the closure of the image of $\bP^2$ under the linear series 
$$
x_1=u_1(u_1+u_2+u_3)^2, \,\, 
x_2=u_2(u_1+u_2+u_3)^2, \,\, 
x_3=u_3(u_1+u_2+u_3)^2,
$$
$$
w=u_1u_2u_3,  
$$
where $\langle u_1,u_2,u_3\rangle =\Gamma(\bP^2,\cO_{\bP^2}(1))$. 
\end{enumerate}

\begin{rem}
There are two isomorphism classes of cubic surfaces with a $\Dfour$ singularity
\cite{bruce-wall} Lemma 4.  The other class is
$$S_0=\{(x_1,x_2,x_3,w): w(x_1+x_2+x_3)^2=x_1x_2(-x_1-x_2) \};$$
it is obtained from $S$ by substituting
$$(w,x_1,x_2,x_3)\mapsto (t^{-2}w,x_1,x_2,
 tx_3+(t-1)x_1+(t-1)x_2)$$
and letting $t\ra 0$ in the resulting equation.

We can distinguish $S$ and $S_0$ geometrically:
In $S$, the three lines not containing $p$ do not share a common point.
In $S_0$, the analogous lines
$$\{w=x_1=0 \}, \{w=x_2=0 \}, \{w=x_1+x_2=0 \} \subset S_0
$$
are coincident at $w=x_1=x_2=0$.
\end{rem}

Let $\beta:\tS\ra S$ denote the minimal desingularization of $S$
and 
$$
\ell_1,\ell_2,\ell_3, m_1, m_2,m_3
$$ 
the strict transforms of the lines. 
The rational map $S\dashrightarrow \bP^2$ induces a morphism 
$\tS\ra \bP^2$
and let $L$ denote the pullback of the hyperplane class. 
Let $E_0,E_1,E_2,E_3$ be the exceptional divisors of $\beta$, ordered 
so that we have the following intersection matrix: 

\begin{equation}
\label{Dmatrix}
\begin{array}{c|ccccccc}
 & L & E_1 & E_2 & E_3 & m_1 & m_2 & m_3 \\
\hline
L   & 1 & 0 & 0 & 0 & 0 & 0 & 0 \\
E_1 & 0 & -2 & 0 & 0  & 1 & 0 & 0   \\
E_2 & 0 & 0 & -2 & 0 & 0 & 1 & 0  \\
E_3 & 0 & 0 & 0 & -2  & 0 & 0  & 1   \\
m_1 & 0 & 1 & 0 & 0 & -1 & 0 & 0   \\
m_2 & 0 & 0 & 1 & 0 & 0 & -1 & 0    \\
m_3 & 0 & 0 & 0 & 1 & 0 & 0 & -1  
\end{array}.
\end{equation}
This is a rank seven unimodular matrix;  since the Picard group 
of $\tS$ has rank seven, it is generated by 
$L,E_1,E_2,E_3,m_1,m_2,m_3$. 
In particular, we have
$$
E_0=L-(E_1+E_2+E_3 +m_1 + m_2+m_3) \,\,\, \text{ and } \,\,\,  
\ell_j=L-E_j-2m_j. 
$$
The anticanonical class is given by 
$$
-K_{\tS}= 3L -(E_1+E_2+E_3) -2(m_1+m_2+m_3) = \ell_1+\ell_2+\ell_3.
$$

\begin{prop}
The effective cone $\NE(\tS)$ is generated by 
$$
\Xi:=\{ E_0, E_1, E_2, E_3, m_j, \ell_j\}.   
$$
\end{prop}
 
\begin{proof}
Each effective divisor $D$ can be expressed as a sum 
$$
D=M_{\Xi}+b_{E_0} E_0 + b_{E_1}E_1 + \ldots + b_{\ell_{3}}\ell_3,
$$
with nonnegative coefficients, where $M_{\Xi}$ intersects each of the
elements in $\Xi$ nonnegatively and thus is in the dual cone
to $\Cone(\Xi)$.  Direct computation shows that
the dual to $\Cone(\Xi)$ has generators
$$
L,  L-E_i-m_i, 2L-E_i -2m_i, 2L-E_i- E_j -2m_i-2m_j,
$$
$$ 
2L-E_i-E_j-m_i-2m_j.
$$
Each of these is contained in $\Cone(\Xi)$:
$$
\begin{array}{rcl}
L                   & =&\ell_i+E_i+2m_i,\\
2L-E_i-2m_i         & =&2\ell_i+E_i+2m_i,\\
2L-E_i-E_j-m_i-2m_j & =&\ell_i+\ell_j+m_i,\\
L-E_i-m_i           &= &\ell_i+m_i, \\
2L-E_i-E_j-2m_i-2m_j&= &\ell_i+\ell_j.
\end{array}
$$
It follows that $M_{\Xi}$ and $D$ are 
sums of elements in $\Xi$ with nonnegative coefficients. 
\end{proof}

Each of the divisors $m_i,\ell_i$ and $E_i$ has a distinguished
nonzero section (up to a constant), denoted $\mu_i, \la_i$ and $\eta_i$, 
respectively.  We have
$$
\{ \la_i\eta_i\mu_i^2, \eta_0\eta_1\eta_2\eta_3\mu_1\mu_2\mu_3\}
\subset \Gamma(L),
$$
and we may identify 
$$
u_i= \la_i\eta_i\mu_i^2\,\,\, \text{  and }\,\,\,  
u_1+u_2+u_3=\eta_0\eta_1\eta_2\eta_3\mu_1\mu_2\mu_3
$$
after suitably normalizing the $\mu_i,\la_i,$ and $\eta_i$.  
The dependence relation among the sections in $\Gamma(L)$ translates 
into 
\begin{equation}
\label{eqn:la}
\la_1\eta_1\mu_1^2+ \la_2\eta_2\mu_2^2+\la_3\eta_3\mu_3^2 = 
\eta_0\eta_1\eta_2\eta_3\mu_1\mu_2\mu_3.
\end{equation}
An argument similar to the one given at the end of Section~\ref{sect:e6}
proves that the natural homomorphism 
$$
k[\eta_0, ..., \eta_3, \mu_i,\la_i ]/  
\langle \la_1\eta_1\mu_1^2+ \la_2\eta_2\mu_2^2+\la_3\eta_3\mu_3^2 - 
\eta_0\eta_1\eta_2\eta_3\mu_1\mu_2\mu_3 \rangle \ra \Cox(\tS)
$$
is an isomorphism. 

\

The cubic surface $S$ admits an $\mathfrak S_3$-action on the coordinates
$x_1,x_2,x_3$. In  particular, it admits nonsplit forms
over nonclosed  ground fields.  
They can be expressed as follows: 
Let $K/k$ be a cubic extension with Galois closure $E/k$.
Fix a basis 
$\{\gamma,\gamma',\gamma''\}$ for $K$ over $k$
so that elements $Y\in K$ can be represented as
$$
Y=y\gamma+ y'\gamma'+y''\gamma''
$$
with $y,y',y''\in k$.  Choose
$\sigma\in \Gal(E/k)$ so that $\sigma$ and $\sigma^2$ are coset
representatives $\Gal(E/k)$ modulo $\Gal(E/K)$.  Then
$$
w\cdot \Tr_{K/k}(Y)^2=\rN_{K/k}(Y)
$$
is isomorphic, over $E$, to $S$:
$$
\begin{array}{ccl}
x_1 & = & Y= y\gamma +y'\gamma'+y''\gamma'' \\
x_2 & = & \sigma(Y)= 
y\sigma(\gamma) +y'\sigma(\gamma')+y''\sigma(\gamma'') \\
x_3 & = & \sigma^2(Y)=
y\sigma^2(\gamma) +y'\sigma^2(\gamma')+y''\sigma^2(\gamma'') 
\end{array}.
$$

Assigning elements $U,V,W\in K$ to 
$\eta_1, \mu_1$ and $\la_1$, respectively, 
the torsor equation \eqref{eqn:la} takes the form
$$
\Tr_{K/k}(UV^2W)=\eta_0\rN_{K/k}(UV).
$$

\bibliographystyle{smfplain}
\bibliography{e6}

\providecommand{\bysame}{\leavevmode ---\ }
\providecommand{\og}{``}
\providecommand{\fg}{''}
\providecommand{\smfandname}{and}
\providecommand{\smfedsname}{\'eds.}
\providecommand{\smfedname}{\'ed.}
\providecommand{\smfmastersthesisname}{M\'emoire}
\providecommand{\smfphdthesisname}{Th\`ese}
\begin{thebibliography}{10}

\bibitem{Ba}
{\scshape V.~V. Batyrev} -- {\og The cone of effective divisors of
  threefolds\fg}, \emph{Proceedings of the International Conference on Algebra,
  Part 3 (Novosibirsk, 1989)}, 
  p.~337--352,
  Contemp. Math., vol. 131, Amer.
  Math. Soc., Providence, RI, 1992.

\bibitem{BP}
{\scshape V.~V. Batyrev {\normalfont \smfandname} O.~N. Popov} -- 
{\og The Cox ring of a Del Pezzo surface \fg}, to appear.  



\bibitem{Br}
{\scshape R.~de~la Bret{\`e}che} -- {\og Nombre de points de hauteur born\'ee
  sur les surfaces de del {P}ezzo de degr\'e 5\fg}, \emph{Duke Math. J.}
  \textbf{113} (2002), no.~3, p.~421--464.


\bibitem{brionproc}
{\scshape M.~Brion {\normalfont \smfandname} C.~Procesi} -- {\og Action d'un
  tore dans une vari\'et\'e projective\fg}, \emph{Operator algebras, unitary
  representations, enveloping algebras, and invariant theory (Paris, 1989),}
  p.~509--539,
  Progr. Math., vol.~92, Birkh\"auser Boston, Boston, MA, 1990.


\bibitem{bruce-wall}
{\scshape J.~W.~Bruce {\normalfont \smfandname} C.~T.~Wall} --
{\og On the classification of cubic surfaces\fg}, 
\emph{J. London Math. Soc. (2)} \textbf{19} (1979), p.~257--267.
 


\bibitem{CLT}
{\scshape A.~Chambert-Loir {\normalfont \smfandname} Y.~Tschinkel} -- {\og On
  the distribution of points of bounded height on equivariant compactifications
  of vector groups\fg}, \emph{Invent. Math.} \textbf{148} (2002), no.~2,
  p.~421--452.

\bibitem{CTS}
{\scshape J.-L. Colliot-Th{\'e}l{\`e}ne {\normalfont \smfandname} J.-J. Sansuc}
  -- {\og La descente sur les vari\'et\'es rationnelles. {II}\fg}, \emph{Duke
  Math. J.} \textbf{54} (1987), no.~2, p.~375--492.

\bibitem{CSS2}
{\scshape J.-L. Colliot-Th{\'e}l{\`e}ne, J.-J. Sansuc {\normalfont \smfandname}
  P.~Swinnerton-Dyer} -- {\og Intersections of two quadrics and {C}h\^atelet
  surfaces. {I}\fg}, \emph{J. Reine Angew. Math.} \textbf{373} (1987),
  p.~37--107.

\bibitem{CSS}
\bysame , {\og Intersections of two quadrics and {C}h\^atelet surfaces.
  {II}\fg}, \emph{J. Reine Angew. Math.} \textbf{374} (1987), p.~72--168.

\bibitem{Cox}
{\scshape D.~A. Cox} -- {\og The homogeneous coordinate ring of a toric
  variety\fg}, \emph{J. Algebraic Geom.} \textbf{4} (1995), no.~1, p.~17--50.

\bibitem{DH}
{\scshape I.~V. Dolgachev {\normalfont \smfandname} Y.~Hu},
	with an appendix by Nicolas Ressayre -- {\og Variation of
  geometric invariant theory quotients\fg}, \emph{Inst. Hautes \'Etudes Sci.
  Publ. Math.} (1998), no.~87, p.~5--56.

\bibitem{fulton}
{\scshape W.~Fulton} -- \emph{Introduction to toric varieties}, Annals of
  Mathematics Studies, vol. 131, Princeton University Press, Princeton, NJ,
  1993.

\bibitem{HT1}
{\scshape B.~Hassett {\normalfont \smfandname} Y.~Tschinkel} -- {\og Geometry
  of equivariant compactifications of {${\bf G}\sb a\sp n$}\fg},
  \emph{Internat. Math. Res. Notices} (1999), no.~22, p.~1211--1230.

\bibitem{HB-C}
{\scshape R.~Heath-Brown} -- {\og 
The density of rational points on Cayley's cubic surface\fg}, 
{\tt math.NT/0210333}, (2003).


\bibitem{HK}
{\scshape Y.~Hu {\normalfont \smfandname} S.~Keel} -- {\og Mori dream spaces
  and {GIT}\fg}, \emph{Michigan Math. J.} \textbf{48} (2000), p.~331--348.

\bibitem{Ko}
{\scshape J.~Koll{\'a}r } -- {\og
Singularities of pairs \fg},
\emph{Algebraic geometry---Santa Cruz 1995}, p.~221--287, 
Proc. Sympos. Pure Math., 62, Part 1, 
Amer. Math. Soc., Providence, RI, 1997. 

\bibitem{KM}
{\scshape J.~Koll{\'a}r {\normalfont \smfandname} S.~Mori},
  with the collaboration of C. H.
  Clemens and A. Corti. -- \emph{Birational
  geometry of algebraic varieties}, Cambridge Tracts in Mathematics, vol. 134,
  Cambridge University Press, Cambridge, 1998.

\bibitem{Mus}
{\scshape M. Musta\c t\u a} --
{\og Vanishing theorems on toric varieties \fg},
\emph{Tohoku Math. J. (2)} \textbf{54} (2002), no. 3, p.~451--470.

\bibitem{Peyre}
{\scshape E.~Peyre} -- {\og Terme principal de la fonction z\^eta des hauteurs
  et torseurs universels\fg}, 
  \emph{Nombre et r\'epartition de points de hauteur born\'ee (Paris, 1996)},
  p.~259--298,
  \emph{Ast\'erisque} (1998), no.~251. 

\bibitem{Peyre2}
\bysame , {\og Torseurs universels et m\'ethode du cercle\fg}, 
  \emph{Rational points
  on algebraic varieties}, 
  p.~221--274,
  Progr. Math., vol. 199, Birkh\"auser, Basel, 2001.

\bibitem{S-ast} 
{\scshape P.~Salberger} -- {\og 
Tamagawa measures on universal torsors and points of bounded
height on Fano varieties\fg}, 
\emph{Nombre et r\'epartition de points de hauteur born\'ee (Paris, 1996)},
p.~91--258,
\emph{Ast\'erisque} (1998), no.~251.

\bibitem{SS}
{\scshape P.~Salberger {\normalfont \smfandname} A.~N. Skorobogatov} -- {\og
  Weak approximation for surfaces defined by two quadratic forms\fg},
  \emph{Duke Math. J.} \textbf{63} (1991), no.~2, p.~517--536.

\bibitem{S}
{\scshape A.~Skorobogatov} -- \emph{Torsors and rational points}, Cambridge
  Tracts in Mathematics, vol. 144, Cambridge University Press, Cambridge, 2001.

\bibitem{Sb}
{\scshape A.~N. Skorobogatov} -- {\og On a theorem of
  {E}nriques-{S}winnerton-{D}yer\fg}, \emph{Ann. Fac. Sci. Toulouse Math. (6)}
  \textbf{2} (1993), no.~3, p.~429--440.

\bibitem{thad}
{\scshape M.~Thaddeus} -- {\og Toric quotients and flips\fg}, \emph{Topology,
  geometry and field theory}, 
  p.~193--213,
 World Sci. Publishing, River Edge, NJ, 1994.

\end{thebibliography}

\end{document}